\documentclass[11pt]{amsart}

\renewcommand{\subjclassname}{\textup{2000} Mathematics Subject
     Classification}

\newtheorem{theorem}{Theorem}[section]
\newtheorem{lemma}[theorem]{Lemma}
\newtheorem{proposition}[theorem]{Proposition}
\newtheorem{corollary}[theorem]{Corollary}

\theoremstyle{definition}

\theoremstyle{remark}
\newtheorem{remark}[theorem]{Remark}

\numberwithin{equation}{section}

\def\R{{\mathbb R}}

\def\C{{\mathbb C}}

\def\Re{{\rm Re}\,}
\def\Im{{\rm Im}\,}
\def\TT{(T(t))_{t\geq 0}}

\def\TcT{({\mathbf T} (t))_{t\geq 0}}

\def\TcT0{({\mathbf T}_0 (t))_{t\geq 0}}

\def\L1{L^1 (\R_+ )}

\newcommand{\XC}{{\mathcal C}}
\newcommand{\XL}{{\mathcal L}}
\newcommand{\XN}{{\mathcal N}}

\newcommand{\e}{{\varepsilon}}

\begin{document}

\title[Optimal polynomial decay]{Optimal polynomial decay of functions and operator semigroups}

\author{Alexander Borichev}
\address{Centre de Math\'ematiques et Informatique,
Universit\'e d'Aix-Marseille I, 39 rue Fr\'ed\'eric Joliot-Curie,
13453 Marseille, France}

\email{borichev@cmi.univ-mrs.fr}

\author{Yuri Tomilov}
\address{Faculty of Mathematics and Computer Science, Nicolas Copernicus University, ul.
Chopina 12/18, 87-100 Toru\' n, Poland and Institute of Mathematics,
Polish Academy of Sciences. \' Sniadeckich str. 8, 00-956 Warsaw,
Poland}

\email{tomilov@mat.uni.torun.pl}
\thanks{The authors were partially supported by the Marie Curie ''Transfer of Knowledge'' programme, project ''TODEQ''.
The first author was also partially supported by the ANR project DYNOP. The second
author was also partially supported by a MNiSzW grant Nr. N201384834.}




\begin{abstract}
We characterize the polynomial decay of orbits of Hilbert space
$C_0$-semigroups in resolvent terms. We also show that results of
the same type for general Banach space semigroups and functions
obtained recently in \cite{BaDu08} are sharp. This settles a
conjecture posed in \cite{BaDu08}.
\end{abstract}

\renewcommand{\subjclassname}{\textup{2000} Mathematics Subject Classification}

\maketitle

\section{Introduction}

One of the main issues in the theory of partial differential equations is to determine
whether the solutions to these equations approach an equilibrium and if yes then
how fast do the solutions approach it. Recently, an essential progress was achieved
in treating such asymptotic problems by operator-theoretical methods involving $C_0$-semigroups.
For different accounts of developments, highlights, and techniques of asymptotic theory of
$C_0$-semigroups see \cite{ABHN01},  \cite{Ba94},  \cite{ChTo07}, and \cite{Ne96}.

In particular, the following result was proved in \cite[p. 803]{Ba90}, see also \cite[p.41-42]{Ba94}.
(The result is implicitly contained already in \cite{ArBa88}.)

\begin{theorem}\label{ingham}
Let $\TT$ be a bounded $C_0$-semigroup on a Banach space $X$ with
generator $A.$ Suppose that $i \mathbb R$ is contained in the
resolvent set $\varrho (A)$ of $A.$ Then
\begin{equation}\label{-11}
\|T(t)A^{-1}\|\to 0, \quad t\to \infty.
\end{equation}
\end{theorem}

In other words, all classical solutions to the abstract Cauchy problem
\begin{equation} \label{cauchy}
\left\{ \begin{array}{ll}
\dot{x} (t) = A x(t) , & t\geq 0 ,\\[2mm]
x(0) = x_0 , & x_0\in X,
\end{array} \right.
\end{equation}
given by $x (t) = T(t)x_0$, $t \ge 0$, $x_0 \in D(A)$, converge
uniformly (on the unit ball of $D(A)$ with the graph norm) to zero
at infinity if $A$ satisfies the conditions of Theorem
\ref{ingham}.

In general, without any additional assumptions, the decay in
\eqref{-11} can be arbitrarily slow. However, in a number of
special situations involving concrete PDE's, e.g. damped
wave equations,  the rate of decay in \eqref{-11} corresponds to the
rate of decay of the energy of the system described by $(T(t))_{t
\ge 0}$, and it is of interest to determine whether this rate of 
decay can be achieved.

 By rewriting equations
in the abstract form \eqref{cauchy}, the rates of the decay of
sufficiently smooth orbits for the corresponding semigroup $\TT$
(and equivalently of solutions to \eqref{cauchy}) can be
associated with the size of the resolvent $R(\lambda,
A)=(\lambda-A)^{-1}$ of $A$ on the imaginary axis. This approach
was initiated in \cite{Le96} and later pursued, in particular,  in
 \cite{Bu98}, \cite{BuHi07}, \cite{Ch09},
 \cite{LeRo97}.
 However, with a few exceptions, the
issue of optimality or (non-optimality) of the rates of growth has
not been studied so far.

The above applications (mainly in the abstract set-up) motivated a
thorough study of the decay rates on $\|T(t)A^{-1}\|$ for bounded
$C_0$-semigroups $\TT$ on Banach spaces in \cite{BaEnPrSchn06},
\cite{LiRa05}, and, most recently, in \cite{BaDu08}. In the latter
paper, the authors developed a unified and simplified approach for
estimating the decay rates for $\|T(t)A^{-1}\|$ in terms of the
growth of $R(is, A)$, $s \in \mathbb R$, using the contour
integrals technique by Newman--Korevaar. In particular, the
following theorem is proved there.

For $A$ as in Theorem \ref{ingham} we define a continuous non-decreasing
function
\begin{equation}\label{mf}
M(\eta)=\max_{t \in [-\eta,\eta]}\|R(it, A)\|,\qquad \eta \ge 0,
\end{equation}
and the associated function
\begin{equation}\label{mflog}
M_{\log}(\eta):=M(\eta)\big(\log(1+M(\eta))+\log (1+\eta)\big),\quad \eta
\ge 0.
\end{equation}
Let $M^{-1}_{\log}$ be the inverse of $M_{\log}$ defined on
$[M_{\log}(0), +\infty).$

\begin{theorem}[Batty, Duyckaerts]\label{battysem}
Let $(T(t))_{t \ge 0}$ be a bounded $C_0$-semi\-group on a Banach
space $X$ with generator $A,$ such that $i \mathbb R \subset
\varrho (A).$ Let the functions $M$ and $M_{\log}$ be defined by
\eqref{mf}
 and \eqref{mflog}. Then
there exist $C,B>0$ such that
\begin{equation}\label{semestgen}
 \|T(t) A^{-1}\|
\le \frac{C}{M^{-1}_{\log}(t/C)}, \quad t \ge B.
\end{equation}
\end{theorem}

Note that in the case $\alpha>0$, $M(\eta) \le C (1+\eta^{\alpha})$, $\eta\ge 0$,
the Batty-Duyckaerts result gives \begin{equation}\label{semest}
 \|T(t) A^{-1}\| \le C
\Bigl( \frac{\log t} {t}  \Bigr)^{\frac{1}{\alpha}}, \quad t \ge B.
\end{equation}
It was conjectured in \cite{BaDu08} that Theorem \ref{battysem} can
be improved by removing the logarithmic factor in \eqref{semest} in
the case when $X$ is a Hilbert space, and that this factor is
necessary if $X$ is merely a Banach space, see \cite[Remark
9]{BaDu08}. (See also the introduction in \cite{BaEnPrSchn06} and
the comments after Theorem 3.5 therein.)

In our paper, we address the problem of {\it optimality} of the
rate of decay in \eqref{semest} and confirm the conjecture from
\cite{BaDu08} for the case of polynomially growing $M$. We show
that the logarithmic factor can be dropped in \eqref{semest} if
$X$ is a Hilbert space. Thus, various results on polynomial decay
of solutions of PDE's, e.g. in \cite{BaEnPrSchn06}, \cite{BuHi07},
 \cite{LiRa05}, \cite{LiRa07}  can be improved to sharp formulations or shown to
be sharp. See also \cite{BaDu08} and references therein.

On the other hand, we show that Theorem \ref{battysem} is
essentially sharp (see Theorems \ref{battysem2} and \ref{xyz}
below). This is done by a function-theoretical construction which
may be interesting for its own sake and may be useful in other
instances related to $C_0$-semigroups as well. We prove, in
particular, that given $\alpha>0$ there exists a Banach space
$X_\alpha$ and a bounded $C_0$-semigroup $(T(t))_{t\ge 0}$ in
$X_\alpha$ with generator $A$ such that
$$
\|R(is,A)\|=O(|s|^\alpha),\qquad |s|\to\infty,
$$
and
$$
\limsup_{t\to\infty}\Bigl(\frac {t}{\log t}\Bigr)^{1/\alpha}\|T(t)A^{-1}\|>0.
$$

The classical problem of estimating the local energy for solutions
to wave equations leads to the study of decay rates for functions
of the form  $\|T_1 T (t) T_2\|,$ where   $T_1, T_2 $ are bounded
operators on  $X,$ so that the assumptions are imposed on the
cut-off resolvent $F(\lambda)= T_1 R(\lambda, A) T_2$  rather than
on the resolvent itself, see e.g \cite{BaDu08}, \cite{Bu98},
\cite{Ch09},
  \cite{Vo04}. Due to the lack of the Neumann series
expansions for $F$ we have to assume that $F$ extends analytically
to the region $\Omega$ of known shape and satisfies certain growth
restrictions there. The domain $\Omega$ and the growth of $F$ in
$\Omega$ are not in general related to each other. Moreover, the
operator-theoretical approach can hardly be used to deal with $F$
since the resolvent identity is not available as well. Thus, it is
natural to put the problem into the framework of decay estimates
for $L^{\infty}(\mathbb R_+, X)$ functions given that their
Laplace transforms extend to the specific $\Omega$ with, in our
case, polynomial estimates. In this direction, we obtain a result
on the polynomial rates of decay for bounded  functions which
partially generalizes \cite[Theorem 10]{BaDu08}. The result has
its version for the rates of decay of $\|T_1 T (t)T_2\|$ thus
improving \cite[Corollary 11]{BaDu08}.

We use standard notations. Given a closed linear operator $A$ we denote
by $\sigma(A)$, $\rho (A)$, ${D}(A)$, and ${\rm Im}\,(A)$ the spectrum
of $A$, the resolvent set of $A$, the domain of $A$ and the image of
$A$ respectively. By $C$, $C_1$ etc. we denote generic constants which may
change from line to line.

The plan of the paper is as follows. In the section \ref{Hilbert} we characterize the rate
of polynomial decay of the semigroup orbits in the Hilbert space via the growth rate
of the resolvent on the imaginary axis. Examples of functions constructed in the section
\ref{functions} show that the version of Theorem \ref{battysem} for functions given in
\cite{BaDu08} is sharp. These examples are used in the section \ref{Banachh} to
show that Theorem \ref{battysem} is itself sharp.
\medskip

\noindent{\bf Acknowledgment}  The authors are grateful to the referees for 
useful comments and suggestions which led to an improvement of the paper.

\section{Decay of Hilbert space semigroups}\label{Hilbert}

We start with recalling two simple and essentially known
statements on $C_0$-semigroups. The first one is a version of the
well-known criterion on the generators of bounded Hilbert space
$C_0$-semigroups, \cite{Go99}, \cite{ShFe00}. In our case, the
proof is particularly easy and, to make our presentation
self-contained, we give an easy argument. Let $\C_+:=\{z\in\C:\Re
z>0 \}$.

\begin{lemma}\label{bound}
Let $(T(t))_{t \ge 0}$ be a $C_0$-semigroup on a Hilbert space $H$
with generator $A.$ Then $(T(t))_{t \ge 0}$ is bounded if and only
if
\begin{eqnarray}
\nonumber & & \C_+ \subset \varrho (A), \mbox{ and} \\
\nonumber & &
\sup_{\xi > 0} \, \xi \int_{\mathbb R} \Bigl(\|R(\xi + i\eta, A)x\|^2 +
\|R(\xi+ i\eta, A^*)x\|^2\Bigr) \; d\eta < \infty \\
\nonumber && \, \mbox{for every } \,  x\in H.
\end{eqnarray}
\end{lemma}

\begin{remark} By the closed graph theorem, in the conditions of Lemma~\ref{bound}, for $x\in H$ we have
$$
\sup_{\xi > 0} \, \xi \int_{\mathbb R} \Bigl( \|R(\xi + i\eta, A)x\|^2 +
\|R(\xi+ i\eta, A^*)x\|^2 \Bigr)\; d\eta \le C\|x\|^2.
$$
\end{remark}

\begin{proof}[Proof of Lemma~{\rm \ref{bound}}]
The necessity is a direct consequence of the Plancherel theorem
applied to the families $\{e^{-\xi t}T(t)x : x \in H, \, \xi > 0 \}$, $\{e^{-\xi t}T^*(t)x : x \in H, \,
\xi > 0 \}$. The sufficiency follows from the
representation
$$
\langle T(t)x, x^*\rangle = \frac{1}{2 \pi i\,t} \int_{\frac 1t - i
\infty}^{\frac 1t + i\infty} e^{\lambda t} {\langle R^2(\lambda, A)
x,x^*\rangle} \, d\lambda, \qquad t>0,
$$
where the integral converges absolutely by the H\"older inequality
and our assumptions, see e.g.  \cite{Go99}, \cite{ShFe00}, \cite{ChTo03}.
\end{proof}

The second  statement allows us to cancel the growth of the resolvent
by restricting it to sufficiently smooth elements of $H$.

\begin{lemma}\label{est}
Let $(T(t))_{t \ge 0}$ be a bounded $C_0$-semigroup on a Hilbert
space $H$ with generator $A,$ such that $i \mathbb R \subset
\varrho (A).$ Then for a fixed $\alpha >0$
$$
\|R(\lambda, A) (-A)^{-\alpha}\| \le C, \qquad \Re \lambda >0,
$$
if and only if
\begin{equation}\label{alpha}
 \|R( is , A) \|= {\rm O}
(|s|^{\alpha}), \qquad s \to \infty.
\end{equation}
\end{lemma}

\begin{proof}
The lemma is proved in  \cite[Lemma 3.2]{LaSh01}, \cite[Lemma
1.1]{HuNe99} in a version saying, in particular, that the condition
$$
\|R(\lambda, A) \| \le C\big(1+ |\lambda|^{\alpha} \big), \qquad 0< \Re \lambda<1,
$$
is equivalent to
$$
\|R(\lambda, A) (-A)^{-\alpha}\| \le C_1, \qquad 0<\Re \lambda<1.
$$
To get our version of the assertion it suffices to apply the maximum principle to the function
$F(\lambda)= R(\lambda,A)\lambda^{-\alpha}\big(1-\frac{\lambda^2}{B^2}\big)$ on the domain
$\{\lambda \in \mathbb C: \Re \lambda \ge 0, \, 1 \le|\lambda | \le B \}$
for large $B$, and to use the estimate
$$
\|R(\lambda, A)\| \le \frac{C}{\Re  \lambda}, \quad \Re\lambda > 0,
$$
for $\lambda$ with $|\lambda|=B$. Then \eqref{alpha} implies
$$
\|R(\lambda, A)\|\le C (1+ |\lambda|^{\alpha}), \quad \Re \lambda > 0.
$$
Since $R(\lambda,A)$ is bounded in any halfplane strictly included in
$\mathbb C_+$, Lemma \ref{est} is reduced to its strip version mentioned above.
\end{proof}

The next theorem is one of the main results of the paper. Its proof is based on a trick: we pass
to a matrix semigroup whose boundedness gives the required rate of decay of $\|T(t)(-A)^{-\alpha}\|$
(and of $\|T(t)A^{-1}\|$) ``for free".

\begin{theorem}\label{main1}
Let  $(T(t))_{t \ge 0}$ be a bounded $C_0$-semigroup on a Hilbert
space $H$ with generator $A$ such that $i \mathbb R \subset
\varrho (A).$ Then for a fixed $\alpha > 0$ the following
conditions are equivalent:
\begin{enumerate}
\item[(i)]
\begin{eqnarray}\label{alpha1}
\|R( is , A) \|= {\rm O} (|s|^{\alpha}), \qquad s \to \infty.
\end{eqnarray}
\item [(ii)]
\begin{eqnarray}\label{decay}
\|T(t)(-A) ^{-\alpha}\|= {\rm O} (t^{-1}), \qquad t \to \infty.
\end{eqnarray}
\item [(iii)]
$$
\|T(t)(-A)^{-\alpha} x \| ={\rm o }(t^{-1}), \qquad \quad t \to
\infty, \, x \in H.
$$
\item [(iv)]
\begin{eqnarray}\label{-1}
\|T(t)A ^{-1}\|= {\rm O} (t^{-1/\alpha}), \qquad t \to \infty.
\end{eqnarray}
\item [(v)]
\begin{eqnarray}\label{-01}
\|T(t)A ^{-1}x\|= {\rm o} (t^{-1/\alpha}), \qquad t \to \infty, x
\in H.
\end{eqnarray}
\end{enumerate}
\end{theorem}

\begin{proof}
The implication ${\rm (ii)}\Rightarrow{\rm (i)}$ was proved in \cite{BaDu08}; 
the implication ${\rm (iii)}\Rightarrow{\rm (ii)}$ is a
consequence of the uniform boundedness principle. Moreover, the
equivalence ${\rm (iv)}\Longleftrightarrow{\rm (ii)}$ was obtained in
\cite[Proposition 3.1]{BaEnPrSchn06} as a consequence of the moment
inequalities for $A$. Its `o'-counterpart ${\rm (iii)}\Longleftrightarrow{\rm (v)}$ can
be obtained by the same argument. Thus, it remains to prove that
${\rm (i)}\Rightarrow{\rm (iii)}$.

${\rm (i)}\Rightarrow{\rm (iii)}$: Let $\mathcal H= H \oplus H$ be the direct sum of
two copies of $H$. Consider the operator ${\mathcal A}$ on $\mathcal H$ given by
the operator matrix
\begin{equation*}
 {\mathcal A}=\left( \begin{array}{cc}
       A  & (-A)^{-\alpha} \\
       O & A \end{array}
       \right)
 \end{equation*}
 with the diagonal domain $ D(\mathcal A)= D(A) \oplus D(A).$
Then $\sigma(A)=\sigma(\mathcal A)$, and the resolvent
$R(\lambda, \mathcal A)$ of  $\mathcal A$ is of the form
\begin{equation*}\label{m}
 R(\lambda, {\mathcal A})=\left( \begin{array}{cc}
       R(\lambda, A)  &   R^2(\lambda, A)(-A)^{-\alpha} \\
             O &  R(\lambda, A) \end{array}
       \right), \qquad \lambda \in \rho (A).
 \end{equation*}
The operator $\mathcal A$ is the generator of the $C_0$-semigroup
$({\mathcal
 T}(t))_{t \ge 0}$ on $\mathcal H$ defined by
 \begin{equation}\label{matr}
 {\mathcal T}(t)=\left( \begin{array}{cc}
       T(t)  &  t T(t)(-A)^{-\alpha} \\
       O &  T(t) \end{array}
       \right),
 \end{equation}
because the resolvents of $\mathcal A$ and of the generator of
$({\mathcal T}(t))_{t \ge 0}$ coincide. By \eqref{alpha1} and Lemma \ref{est},
$$
\|R(\lambda, A)(-A)^{-\alpha}\| \le C, \quad \Re \lambda >0.
$$
Hence, for every $\mathbf x=(x_1, x_2) \in \mathcal H$ and $\lambda \in \mathbb C_+$,
\begin{equation}\label{maj}
\| R(\lambda, {\mathcal A}){\mathbf x}  \|^2 \le C \left(\|
R(\lambda, A)x_1 \|^2 + \| R(\lambda, A) x_2 \|^2 \right),
\end{equation}
and similarly
\begin{equation}\label{maj*}
\| R(\lambda, {\mathcal A^*}){\mathbf x}  \|^2 \le C \left(\|
R(\lambda,  A^*)x_1 \|^2 + \| R(\lambda,  A^*) x_2 \|^2 \right).
\end{equation}
By Lemma \ref{bound},
$$
\sup_{\xi > 0} \, \xi \int_{\mathbb R} \Bigl(\|R(\xi + i\eta, A)x\|^2 + \|R(\xi+ i\eta, A^*)x\|^2\Bigr) \; d\eta < \infty
$$
for every $x\in H$, so that
$$
\sup_{\xi > 0}\,  \xi \int_{\mathbb R} \Bigl(\|R(\xi + i\eta, \mathcal A) \mathbf x\|^2 +
\|R(\xi+ i\eta, \mathcal A^*) \mathbf x\|^2 \Bigr) \; d\eta < \infty
$$
for every $\mathbf x\in \mathcal H$. Then again by Lemma \ref{bound},
$({\mathcal T}(t))_{t \ge 0}$ is bounded on $\mathcal H$.
Since $\TT$ is bounded, the definition of $({\mathcal T}(t))_{t \ge 0}$
implies that
$$
\sup_{t \ge 0} \|t T(t) (-A)^{-\alpha}\| < \infty.
$$
Furthermore, $i\mathbb R\subset \rho (\mathcal A)$.  Then by
Theorem \ref{ingham},
\begin{equation}\label{stab}
{\mathcal T}(t)x \to 0, \qquad t\to \infty, \quad \text {for every }x \in H,
\end{equation}
since ${D}(\mathcal A)={\rm Im}\,({\mathcal A^{-1}})$ is dense in
$\mathcal H$.
Again by Theorem \ref{ingham}, \eqref{stab} implies that
$$
\|t T(t) (-A)^{-\alpha}x\|={\rm o}(1), \qquad t \to \infty,\quad x \in H.
$$
\end{proof}

\begin{remark}
If one is merely interested in the proof of the implication  ${\rm
(i)}\Rightarrow{\rm (ii)}$, then there is no need to invoke Theorem
\ref{ingham} as the argument above shows.

There is another argument for ${\rm (i)}\Rightarrow{\rm (iii)}$,
which does not use Theorem \ref{ingham} too. Let $n\ge 1+\alpha$
be an integer. By the resolvent identity, $R(\cdot,A)x \in L^2
(i\mathbb R,H)$, and, moreover, $R(\cdot,A)x \in H^2 (\mathbb
C_+,H)$ (the Hardy class in the right half-plane) for every $x$
from the dense set ${\rm Im}\,(-A)^{-n}$. Then,
\begin{equation}\label{limit}
\lim_{\xi \to 0+} \xi \int_{\mathbb R} \|R(\xi + i\eta,A)x\|^2 \; d\eta=0, \qquad x \in
{\rm Im}\,(-A)^{-n},
\end{equation}
and since $(T(t))_{t \ge 0}$ is bounded, the last relation holds
for all $x \in H$. 
By the simple integral resolvent stability criterion
from \cite[Theorem 3.1]{To01},
this means that $(T(t))_{t \ge 0}$ satisfies
\eqref{stab}.
 Using the estimates \eqref{maj}, \eqref{maj*}, we conclude that
 $R(\xi + i\eta, \mathcal A)x$, $\xi >0$, satisfies an analog of \eqref{limit} for
every $x \in \mathcal H$. Since the semigroup  $({\mathcal
T}(t))_{t \ge 0}$ is bounded, by the same stability criterion,
$({\mathcal T}(t))_{t \ge 0}$ is stable. Therefore, in particular,
$$
\|T(t)(-A)^{-\alpha} x \| ={\rm o }(\frac{1}{t}), \qquad  t \to
\infty, \, x \in H.
$$
\end{remark}

\begin{remark}\label{normal}
Consider the elementary example of a multiplication
$C_0$-semi\-group  $(T(t))_{t \ge 0}$,
$$
(T(t)f)(z)=e^{tz}f(z), \qquad t \ge 0,
$$
on $L^2 (S, \mu),$ where $S := \{z \in \mathbb C: \Re z < -1/(1+|\Im z|)^\alpha \}$ and $\mu$ is
Lebesgue measure on $S$.

The operator $A f(z)=z f(z)$ with  maximal domain is the generator of $(T(t))_{t \ge 0}$, and
for $f,g\in L^2(S, \mu)$ we have
\begin{eqnarray*}
\langle R(\lambda, A)f,g\rangle &=&
\int_S \frac{f(\zeta)\overline{g(\zeta)}\,d\mu(\zeta)}{\lambda -\zeta}, \qquad \lambda \in \mathbb C \setminus S,\\
\langle T(t)(-A)^{-\alpha}f,g\rangle &=& \int_S \frac{e^{t\zeta}f(\zeta)\overline{g(\zeta)}\,d\mu(\zeta)}
{(-\zeta)^\alpha}, \qquad t \ge 0.
\end{eqnarray*}
Straightforward estimates give that
$$
\|R(is, A)\|={\rm O}(|s|^{\alpha}), \,|s| \to \infty, \quad \text{and}
\qquad \lim_{t\to +\infty}t\|T(t)(-A)^{-\alpha}\|=1,
$$
which demonstrates the optimality of the rates of decay 
in Theorem \ref{main1}. A similar fact is
true for the decay rate of $\|T(t)A^{-1}x\|$, $x \in H$, see, for
example, \cite[Proposition 4.1]{BaEnPrSchn06}.
\end{remark}

\begin{remark}
An advantage of the construction in the proof of Theorem
\ref{main1} is that it reduces the problem of finding optimal
polynomial rates in \eqref{-1} under resolvent growth conditions  to
proving the boundedness of  $({\mathcal T}(t))_{t \ge 0}$ under the same
kind of conditions. While there is a criterion for boundedness of
Hilbert space semigroups in terms of boundary behavior of certain
resolvent means, see e.g. \cite{Go99}, \cite{ShFe00}, a similar
criterion for semigroups on Banach spaces is yet to be found. Note
that the Banach space semigroups $(T(t))_{t \ge 0}$ satisfying
$\|T(t)\| \le w(t),$ where $w:\mathbb R_+ \to \mathbb R_+$ is
submultiplicative, can be {\it characterized} by
Hille-Yosida type conditions, see \cite[Theorem 5.1]{Cho02} and the
comments following it. However, we do not yet know how to verify such conditions
for $({\mathcal T}(t))_{t \ge 0}$ and for the weight $w(t)=C \ln (e+t)$ to conclude that
$\|{\mathcal T}(t)\| \le C \ln (e+t)$, $t\ge 0$, recovering thus Theorem \ref{battysem} for
polynomially growing functions $M$ in the Banach space setting.
\end{remark}

\section{Decay of functions}\label{functions}

Let $M$ be a continuous non-decreasing function, and let $M_{\rm log}$ be defined by \eqref{mflog}.

The following result is an analog of Theorem \ref{battysem} for the decay
of functions $f \in L^{\infty}(\mathbb R_+,X)$ (see also \cite[Theorem 10]{BaDu08}).
Given $f\in L^{\infty}(\mathbb R_+, X)$, its Laplace transform is defined by
$$
\widehat f (z)=\int_0^{\infty}e^{-zt}f(t)\,dt.
$$


\begin{theorem}[Batty, Duyckaerts]\label{battysem2}
Let  $f \in L^{\infty}(\mathbb R_+, X)$ be
such that
\begin{enumerate}
\item [a)]$ \widehat f$ extends analytically to the domain
\newline\noindent
$\Omega:=\{z\in\mathbb C: \Re z > -1/M(|\Im z|)\}$, and
\item[b)]
$$
\|\widehat f (z)\| \le  M(|\Im z|), \qquad z\in \Omega.
$$
Then there exist $C,C_1>0,t_0$ {\rm(}depending on $\|f\|$ and $M${\rm)} such that
$$
\Bigl\|\widehat f(0) -\int_{0}^{t} f(s) \, ds \Bigr\| \le
\frac{C_1}{M^{-1}_{\log}(t/C)}, \qquad t\ge t_0.
$$
\end{enumerate}
\end{theorem}

If we are interested only in polynomial rate of growth of
$\widehat f$ and in (possibly different) polynomial rate of
narrowing of $\Omega\setminus \mathbb C_+$, we can formulate the
following variant of  \cite[Theorem 10]{BaDu08}.

\begin{proposition}\label{x1}
Let $\alpha,\beta,c_1,c_2>0$, $\Omega:=\{z\in\mathbb C: \Re z
> -c_1(1+|\Im z|)^{-\alpha}\}$. Suppose that $f \in
L^{\infty}(\mathbb R_+, X)$
is such that $\widehat f$ admits an analytic extension to $\Omega$
and
$$
\|\widehat f (z)\| \le
c_2(1+|\Im z|)^\beta, \qquad z\in \Omega.
$$
Then
$$
\Bigl\|\widehat f(0) -\int_{0}^{t} f(s) \, ds \Bigr\| \le
C\Bigl(\frac{\log t}{t}\Bigr)^{1/\alpha}, \qquad t\ge 2,
$$
with $C$ depending only on $c_1, c_2, \|f\|_{\infty},\alpha, \beta$.
\end{proposition}

\begin{remark}
Thus, given polynomial growth of $\widehat  f$ in $\Omega$, the rate of
polynomial decay of  $(\widehat f(0) -\int_{0}^{t} f(s) \, ds)$ is determined only by the shape of
$\Omega.$ This is in sharp contrast to the semigroup case where the
growth of $R(\lambda, A)$ in $\Omega$ and the size of $\Omega$ are
related to each other due to the Neumann series expansions of the
resolvent.
\end{remark}

Without loss of generality we assume that $\widehat f$ is continuous up to $\partial\Omega$.

\begin{lemma}
\label{unn} 
Under the conditions of Proposition~{\rm \ref{x1}}, for any
$\varepsilon>0$ and  for some $C_1,C_2$ depending only on $c_1,
c_2, \|f\|_{\infty},\alpha, \beta,\varepsilon$ we have
$$
\|\widehat f (z)\| \le C_2(1+|\Im z|)^{\alpha+\varepsilon}, \qquad
z\in \Omega',
$$
where $\Omega':=\{z\in\mathbb C: \Re z> -C_1(1+|\Im z|)^{-\alpha}\}$.
\end{lemma}

\begin{proof} Suppose that $\beta>\alpha+\varepsilon$ (otherwise, there is nothing to prove).  We use the fact that the function $\log
\|\widehat f\|$ is subharmonic in $\Omega$. Fix
$A>2\beta/\varepsilon$ and $y>1$ (we deal with the case $y<-1$ by
symmetry), and set
\begin{gather*}
Q_y=\{z:-c_1(y+1)^{-\alpha}<\Re z<\frac
{c_1}{A}(y+1)^{-\alpha},\,
y-1<\Im z<y+1\},\\
E_1=\partial Q_y\cap\Bigl(\frac {c_1}{A}(y+1)^{-\alpha}+i\mathbb R\Bigr),\\
E_2=\partial Q_y\setminus E_1.
\end{gather*}
Next we use that
\begin{gather*}
\log\|\widehat f(z)\|\le \log c_2+\beta\log (y+1),\qquad z\in E_2,\\
\log\|\widehat f(z)\|\le \log \frac
{\|f\|_{\infty}\cdot A}{c_1}+\alpha\log (y+1),\quad z\in E_1,
\end{gather*}
and an estimate of harmonic measure in the thin rectangle $Q_y$,
\begin{equation}
\omega\bigl(-\frac {c_1}{A}(y+1)^{-\alpha}+iy,E_2,Q_y\bigr)<\frac 2A,
\label{zZ}
\end{equation}
for large $y$; here $\omega(z,E,U)$ is harmonic measure of
$E\subset \partial U$ with respect to $z\in U$. By the theorem on
two constants (see, for example, \cite[VII B1]{Koo}) we obtain
\begin{gather*}
\log\|\widehat f(-\frac {c_1}{A}(y+1)^{-\alpha}+iy)\|\le
\omega(-\frac {c_1}{A}(y+1)^{-\alpha}+iy,E_2,Q_y)
\sup_{E_2}\log\|\widehat f\|\\+ \omega(-\frac
{c_1}{A}(y+1)^{-\alpha}+iy,E_1,Q_y) \sup_{E_1}\log\|\widehat f\|.
\end{gather*}
Hence there exists $y_0$ (depending on $c_1,c_2, \alpha, \beta,
A,$ and on $\|f\|_\infty$) such that
$$
\log\|\widehat f(-\frac {c_1}{A}(y+1)^{-\alpha}+iy)\|\le
(\alpha+\varepsilon)\log(y+1), \quad y > y_0.
$$
To get \eqref{zZ} we could just observe that on $Q_y$ we have
$$
\omega(\cdot,E_1,Q_y)\ge \omega(\cdot,E_1,Q)-\omega(\cdot,E',Q),
$$
where $Q=\{z:-c_1(y+1)^{-\alpha}<\Re z<\frac
{c_1}{A}(y+1)^{-\alpha}\}$ and $E'=\partial Q\setminus \partial Q_y$, and use elementary
estimates  for  $\omega(\cdot,E_1,Q)$ and $\omega(\cdot,E',Q)$
obtained via conformal mapping of the strip $Q$ onto the
half-plane.

\end{proof}

\begin{proof}[Proof of Proposition~{\rm \ref{x1}}] By Lemma~\ref{unn}, we can assume that $\beta\le 2\alpha$.

Next, we follow the argument from \cite{BaDu08} (which goes back to
\cite{Ko82} and \cite{Ne80}). By the Cauchy integral formula, for any sufficiently small contour
$\gamma$ around $0$ and for any $R>1$ we have
\begin{equation}
\widehat f(0) -\int_{0}^{t} f(s) \, ds= \frac 1{2\pi i}\int_\gamma
\Bigl(1+\frac {z^2}{R^2}\Bigr) \Bigl(\widehat f(z) -\int_{0}^{t}
e^{-zs}f(s) \, ds\Bigr)e^{zt}\,\frac{dz}{z}.
\label{unn1}
\end{equation}
Let $\gamma_1=R\mathbb T\cap \mathbb C_+$, $\gamma_2=R\mathbb T\cap
(\Omega\setminus \mathbb C_+)$, $\gamma_3=R\mathbb T\setminus
\mathbb C_+$, $\gamma_4=\partial\Omega\cap R\mathbb D$,
where $\mathbb D$ is the unit disc of $\mathbb C$, and $\mathbb T=\partial \mathbb D$. Then
\begin{gather*}
2\pi \Bigl\|\widehat f(0) -\int_{0}^{t} f(s) \, ds \Bigr\| \le
\Bigl\|\int_{\gamma_1} \Bigl(1+\frac {z^2}{R^2}\Bigr) \Bigl(\widehat
f(z) -\int_{0}^{t} e^{-zs}f(s) \, ds\Bigr)
e^{zt}\,\frac{dz}{z}\Bigr\|\\+ \Bigl\|\int_{\gamma_2} \Bigl(1+\frac
{z^2}{R^2}\Bigr) \widehat f(z) e^{zt}\,\frac{dz}{z}\Bigr\|+
\Bigl\|\int_{\gamma_3} \Bigl(1+\frac {z^2}{R^2}\Bigr)
\Bigl(\int_{0}^{t} e^{-zs}f(s) \, ds\Bigr)
e^{zt}\,\frac{dz}{z}\Bigr\|\\+ \Bigl\|\int_{\gamma_4} \Bigl(1+\frac
{z^2}{R^2}\Bigr) \widehat f(z)
e^{zt}\,\frac{dz}{z}\Bigr\|=I_1+I_2+I_3+I_4.
\end{gather*}
Now
\begin{gather*}
I_1\le c\int_{-\pi/2}^{\pi/2} \Bigl(\int_{0}^\infty e^{-sR\cos
\theta}\|f(s+t)\|\,ds\Bigr)\, \cos\theta\,d\theta
\le \frac{c\cdot \| f\|_{L^{\infty}(\mathbb R_+, X)}}{R},\\
I_2\le c\int_{0}^{c_1 (R+1)^{-\alpha}}
\frac{c_2(R+1)^\beta y\,dy}{R^2}
 \le cc^2_1c_2(R+1)^{\beta-2\alpha-2},\\
I_3\le c\int_{-\pi/2}^{\pi/2} \Bigl(\int_{0}^t e^{-sR\cos
\theta}\|f(t-s)\|\,ds\Bigr)\, \cos\theta\,d\theta
\le \frac{c\cdot \|f\|_{L^{\infty}(\mathbb R_+, X)}}{R},\\
I_4\le c\int_{0}^{R}c_2(s+1)^\beta
e^{-c_1t(s+1)^{-\alpha}}\frac{ds}{s+1}.
\end{gather*}
Setting $R=c(c_1,\alpha,\beta)(t/\log t)^{1/\alpha}$, we obtain
$$
I_1+I_2+I_3+I_4\le C(c_1, c_2,  \| f \|_{\infty},
\alpha,\beta)\Bigl(\frac{\log t}{t}\Bigr)^{1/\alpha}, \qquad t\ge 2.
$$
\end{proof}

\begin{remark}
We can avoid using Lemma~\ref{unn} by replacing the term $(1+\frac {z^2}{R^2})$
in the right hand side of \eqref{unn1} by $(1+\frac {z^2}{R^2})^N$ for a suitable $N$.
However, Lemma~\ref{unn} provides some additional information on the growth of $\widehat f$
close to the imaginary axis.
\end{remark}

The statement of Proposition~\ref{x1} can be sharpened if $f$ belongs to the space
of boun\-ded uniformly continuous $X$-valued functions ${\rm
BUC}(\mathbb R_+, X)$.

\begin{proposition}\label{smallo}
Suppose that $f \in {\rm BUC}(\mathbb R_+, X)$ satisfies the
assumptions  of Proposition {\rm \ref{x1}}. Then
\begin{equation}
\Big\|\widehat f(0) -\int_{0}^{t} f(s) \, ds \Big\|^{\alpha} = {\rm
o} \Bigl(\frac{\log t}{ t} \Bigr), \quad t \to \infty.
\label{dop1}
\end{equation}
\end{proposition}

\begin{proof}
By \cite[Corollary 4.4.6]{ABHN01}, our hypothesis on $f$ imply
that
 \begin{equation}\label{o}
\| f(t)\|={\rm o}(1), \quad t\to\infty.
\end{equation}

Arguing as in the proof of Proposition~\ref{x1}, and using
\eqref{o} we obtain
\begin{gather*}
I_1\le \frac{c}{R}\sup_{s\ge 0}\|f(s+t)\|=o\Bigl(\frac{1}{R}\Bigr),\qquad t\to\infty,\\
I_3\le c\sup_{-\pi/2<\theta<\pi/2}\cos\theta\, \int_{0}^t e^{-sR\cos
\theta}\|f(t-s)\|\,ds=o\Bigl(\frac{1}{R}\Bigr),\qquad t\to\infty,
\end{gather*}
and \eqref{dop1} follows.
\end{proof}

The following result is a version of Proposition \ref{x1} written in
the semigroup language. Let $\alpha,\beta>0$ and the domain $\Omega$ be as above.

\begin{corollary}\label{semigroup1}
Let $(T(t))_{t \ge 0}$ be a $C_0$-semigroup on a Banach space $X$
with  generator $A,$ such that $\sup_{t \ge 0}\|T(t)\|:= M <
\infty$, and let $T_1$ and $T_2$ be bounded operators on $X$.
Suppose that $F(z)=T_1 R(z, A)T_2$ admits a holomorphic
extension to the domain $\Omega$,
and suppose that
$$
\|F(z)\|\le M_1(1+|\Im z|)^\beta, \qquad z\in \Omega.
$$
Then there exist $c=c (M, M_1, \alpha,\beta)>0$ such that
\begin{equation}
\|T_1 T(t)(I-A)^{-1}T_2\| \le c \Bigl(\frac{\log
t}{t}\Bigr)^{1/\alpha}, \qquad t\ge 2.
\end{equation}
\end{corollary}

\begin{proof}
As in the proof of \cite[Corollary 11]{BaDu08} we consider
the function
\begin{eqnarray*}
f(t)&=&\frac{d}{dt} \bigl( T_1 T(t)(I-A)^{-1}T_2\bigr)= T_1
T(t)A(I-A)^{-1}T_2\\
&=& - T_1 T(t)T_2 + T_1 T(t)(I-A)^{-1}T_2.
\end{eqnarray*}
Its Laplace transform $\widehat f$ extends analytically to $\Omega$ with the
same estimates as $F$ (up to a constant multiple). Since
$$
\widehat f(0) -\int_{0}^{t} f(s) \, ds=-T_1 T(t)(I-A)^{-1}T_2,
$$
Proposition~\ref{x1} yields the claim.
\end{proof}

Next result shows that Proposition~\ref{smallo} is sharp,
at least for $\beta>\alpha/2$. It suffices to consider just
scalar-valued functions.

\begin{theorem}\label{example}
Given $\alpha>0$, $\beta>\alpha/2$, and a positive function $\gamma
\in C_0(\mathbb R_+)$, there exists a function $f\in C_0 (\mathbb
R_+)$ such that
\begin{enumerate}
\item [a)]$ \widehat f$ admits an analytic extension to the region
$\Omega:=\{z\in\mathbb C: \Re z >  -1/(1+|\Im z|)^\alpha\}$
and
$$
\widehat f (z) (1 +  |\Im z|)^{-\beta} \to 0, \quad |z| \to
\infty, \,z\in \Omega,
$$
and
\item [b)]
$$
\limsup_{t \to \infty} \frac{t}{\gamma(t) \log t} \Big|\widehat f(0)
-\int_{0}^{t} f(s) \, ds \Big|^{\alpha} > 0.
$$
\end{enumerate}
\end{theorem}
To prove Theorem \ref{example} we imitate the multiplication
semigroup example described in Remark \ref{normal}. However, we
choose an infinite charge $\mu$ so that while formal integral
expressions for the resolvent, the semigroup and its orbits remain
the same as in Remark \ref{normal}, the corresponding size estimates 
behave in a different way due to the lack of absolute
convergence of the integrals.

The proof of Theorem \ref{example} is based on the following
lemma. Denote by $\chi_A$ the characteristic function of a set $A
\subset \mathbb C.$
\begin{lemma}\label{prop} Given $Q>0$ and $\e>0$, there exist
an integer $k>Q$ and a complex measure $\mu$ with compact
support in $\mathbb C\setminus \Omega$ such that for some
$B=B(\alpha,\beta)$, we have
\begin{eqnarray}
\Bigl|\int_{\mathbb C \setminus
\Omega}\frac{d\mu(\zeta)}{z-\zeta}\Bigr|&\le&
B(1+|\Im z|)^{\beta}\cdot\chi_{\{z:|z|>Q\}}+\e, \qquad z\in \Omega,
\label{X1}\\
\Bigl|\int_{\mathbb C \setminus \Omega}
e^{t\zeta}\,d\mu(\zeta)\Bigr|&\le& B\chi_{\{t:t>Q\}}+\e,
\qquad t\ge 0,\label{X3}\\
1/B&\le& \Bigl|\int_{\mathbb C \setminus \Omega} e^{k\zeta}\,d\mu(\zeta)\Bigr|\le B,\label{X4}\\
\Bigl|\int_{\mathbb C \setminus \Omega} e^{t
\zeta}\,\frac{d\mu(\zeta)}{\zeta}\Bigr|^\alpha &\le& B\frac{\log
t}{t}\cdot\chi_{\{t:Q<t<2k\}}+\frac{\e}{t+1},\qquad t\ge 0,\label{X5}
\end{eqnarray}
and
\begin{equation}
\Bigl|\int_{\mathbb C \setminus \Omega} e^{k
\zeta}\,\frac{d\mu(\zeta)}{\zeta}\Bigr|^\alpha \ge \frac{\log
k}{Bk}\label{X6}.
\end{equation}
\end{lemma}

\begin{proof}
Choose $H>Q$ large enough, and for an integer $k$, $k \ge 2$,
such that
$$
H^{\alpha}\le k\le H^{3\alpha/2},
$$
define
\begin{eqnarray*}
A&:=&2k\log k, \\
\tau&:=&A^{k-1}/\sqrt k,\\
q&:=&e^{2\pi i/k}, \\
w&:=&iH-1.
\end{eqnarray*}

We define also a finite measure
$$
\mu:=\tau\sum_{1 \le s \le k}q^s(1+q^{s}/(Aw))\delta_{w+q^s/A},
$$
where $\delta_x$ is the unit mass at $x$.

Observe that
\begin{eqnarray}\label{identity}
\sum_{1 \le s \le k} \frac{q^s}{x -q^s}= \frac{k}{ x^k-1}, \qquad
\sum_{1 \le s \le k} \frac{q^{2s}}{x -q^s}&=& \frac{k x }{x^k-1}.
\end{eqnarray}
Indeed, to prove the first identity we use that
$$
\sum_{1 \le s \le k}\frac{q^s}{x-q^s}= \frac {P(x)}{x^k-1},
$$
for some polynomial $P$ with ${\rm deg} \, P < k$ such that
$P(x)=P(qx)$ so that $P(x)={\rm const}$, and then $P(x)=P(0)=k$. The
second equality can be proved in a similar way by using that if
$$
\sum_{1 \le s \le k} \frac{q^{2s}}{x-q^s}=\frac{P(x)}{x^k-1},
$$
then $P(qx)=qP(x)$ which yields  $P(x)=k x$.

Now using \eqref{identity} we get
\begin{eqnarray*}
\XC\mu(z)\overset{\rm def}{=}\int_{\mathbb C \setminus
\Omega} \frac{d\mu(\zeta)}{z-\zeta}&=&
\tau \sum_{1 \le s \le
k}\frac{q^s(1+q^{s}/(Aw))}{(z-w)-q^s/A}\\
&=& \tau \sum_{1 \le s \le k}\Bigl[ \frac{Aq^s}{A(z-w)-q^s}+
\frac{q^{2s}/w}{A(z-w)-q^s} \Bigr]
\\&=&
\tau A\frac{k}{A^k(z-w)^k-1}+\frac\tau w\frac{kA(z-w)}{A^k(z-w)^k-1}
\\&=& \frac z w\frac{kA\tau}{A^k(z-w)^k-1}.
\end{eqnarray*}

If $z\in \Omega$, $H>1$, then
$$
|z-w|\ge 1-\frac 1{H^\alpha}.
$$
Indeed, if $\Im z<H-1$, then $|z-w|>1$, and otherwise, $\Re z\ge -H^{-\alpha}$.

Now if $\frac H2\le \Im z\le |z|\le 2 H$, $z\in \Omega$, then we use that
$A|z-w|>2$, $|A^k(z-w)^k-1|>A^ke^{-k/H^{\alpha}}/2$, to obtain that
$$
\Bigl|\frac z w\frac{kA\tau}{A^k(z-w)^k-1}\Bigr|\le c \tau kA^{1-k}
e^{k/H^{\alpha}}=c\sqrt k e^{k/H^{\alpha}}.
$$
From now on we assume that $k$ satisfies the condition
\begin{equation}
\sqrt k e^{k/H^{\alpha}}\le H^{\beta}. \label{1}
\end{equation}
Then
\begin{equation}
|\XC\mu(z)|\le cH^{\beta}, \qquad \frac H2\le \Im z\le |z|\le
2H, \, \, z\in \Omega. \label{X2}
\end{equation}

If $z\in \Omega$ and $|\Im z-H|+\bigl||z|-H\bigr|>H/2$, then
$|z-w|>c\max(|z|,H)$, and under condition \eqref{1} we have for
large $H$:
$$
|\XC\mu(z)|\le \frac{c_1|z| \sqrt k}{H ( c\max(|z|,H))^k} \le
\frac{c_1|z| \sqrt k}{H( c\max(|z|,H))^{\alpha+1}} \le \e.
$$
This, together with \eqref{X2}, proves \eqref{X1}.

Next,
$$
\XL\mu(t)\overset{\rm def}{=} \int_{\mathbb C \setminus
\Omega} e^{t\zeta}\,d\mu(\zeta)=
\tau\sum_{1 \le s \le k}q^s(1+q^{s}/(Aw))e^{t(w+q^s/A)},
$$
and
\begin{equation}
|\XL\mu(t)|=\tau e^{-t} \Bigl|\sum_{1 \le s \le
k}q^s(1+q^{s}/(Aw))e^{q^st/A}\Bigr|.
\label{dop2}
\end{equation}
Furthermore, we have 
\begin{eqnarray*}
&&\sum_{1 \le s \le k}q^s(1+q^{s}/(Aw))e^{q^st/A}\\&=& \sum_{1 \le s
\le k}\sum_{n\ge 0}(q^s+q^{2s}/(Aw))(q^st/A)^n\frac 1{n!}
\\&=&
k \sum_{m\ge 1}\Bigl[\frac{t^{km-1}}{A^{km-1}}\cdot \frac
1{(km-1)!}+ \frac{t^{km-2}}{A^{km-2}}\cdot \frac 1{(km-2)!}\cdot
\frac 1{Aw} \Bigr]\\&=& \frac{kt^{k-1}}{A^{k-1}(k-1)!}\sum_{m\ge 1}
\Bigl(\frac{t^k}{A^k}\Bigr)^{m-1}\frac{(k-1)!}{(km-1)!} \Bigl(1+
\frac {km-1}{tw}\Bigr),
\end{eqnarray*}
and
$$
|\XL\mu(t)|= \frac{k^{3/2}t^{k-1}e^{-t}}{k!}\sum_{m\ge 1}
\Bigl(\frac{t^k}{A^k}\Bigr)^{m-1}\frac{(k-1)!}{(km-1)!} \Bigl(1+
\frac {km-1}{tw}\Bigr).
$$

Thus, for some constants $c,c_1,c_2,c_3$ we have
\begin{multline*}
|\XL\mu(t)| \le c  \frac{k^{3/2}t^{k-2}e^{-t}}{k!}\sum_{m\ge 1}
\Bigl(\frac{k}{HA}\Bigr)^{k(m-1)}\frac{(k-1)!}{(km-1)!} \frac {km}{H}\\
\le c_1 e^{-t}k^{5/2}t^{k-2}/(k!H),\qquad  0\le t\le k/H,
\end{multline*}
and
$$
c_2e^{-t}k^{3/2}t^{k-1}/k!\le
|\XL\mu(t)| \le c_3 e^{-t}k^{3/2}t^{k-1}/k!,\qquad k/H\le t\le A.
$$
If $t\ge A$, then by \eqref{dop2},
$$
|\XL\mu(t)| \le 2\tau e^{-t}k e^{t/A}=2 \sqrt{k} A^{k-1} e^{-t(1-1/A)}.
$$

The function $t\mapsto e^{-t}t^{k-2}$ attains its maximum on
$[0,k/H]$ at $t=k/H$; the function $t\mapsto e^{-t}t^{k-1}$
attains its maximum on $[k/H,A]$ at $t=k-1$; the function
$t\mapsto e^{-t(1-1/A)}$ attains its maximum on $[A,+\infty)$ at
$t=A$. Using that $A=2 k\log k$, by Stirling's formula, we obtain
that
\begin{gather*}
0< c_1\le |\XL\mu(k)|\le c_2 |\XL\mu(k-1)|\le c_3,\\
\max_{[k/H,k/2]}|\XL\mu|={\rm o}(1), \qquad H\to\infty,\\
e^{-k/H}k^{5/2}(k/H)^{k-2}/(k!H)\to 0,\qquad H\to\infty,\\
\sqrt{k} A^{k-1} e^{-A}\to 0,\qquad H\to\infty.
\end{gather*}

Hence,
\begin{gather*}
0< c_1\le |\XL\mu(k)|\le c_2 \max_{\mathbb R_+}|\XL\mu|\le c_3,\\
\max_{[0,k/2]\cup[A,+\infty)}|\XL\mu|={\rm o}(1), \qquad H\to\infty,
\end{gather*}
and \eqref{X3} and \eqref{X4} follow for large $H$.

Finally,
$$
\XN\mu(t)\overset{\rm def}{=} \int_{\mathbb C \setminus
\Omega}
e^{t\zeta}\,\frac{d\mu(\zeta)}{\zeta}= \frac{\tau
e^{-t+iHt}}{w}\sum_{1 \le s \le k}q^se^{q^st/A}.
$$
Here we use the formula
\begin{eqnarray*}
\sum_{1 \le s \le k}q^se^{q^st/A}&=& \sum_{1 \le s \le k}\sum_{n\ge
0}q^s(q^st/A)^n\frac 1{n!} \\
&=& k\sum_{m\ge1}\frac{t^{km-1}}{A^{km-1}}\cdot \frac
1{(km-1)!}\\&=& \frac{kt^{k-1}}{A^{k-1}(k-1)!} \sum_{m\ge
1}\Bigl(\frac{t^k}{A^k}\Bigr)^{m-1}\frac{(k-1)!}{(km-1)!}
\\&=&(1+{\rm o}(1))\frac{kt^{k-1}}{A^{k-1}(k-1)!},\qquad 0\le t\le A,\,\,H\to\infty.
\end{eqnarray*}
Therefore, by Stirling's formula,  we have for large $H$:
\begin{equation}
|\XN\mu(k)|\ge c_1\frac{\tau e^{-k}}{H}\cdot
\frac{k^{k}}{A^{k-1}(k-1)!} \ge \frac {c_2}H. \label{X8}
\end{equation}
Moreover,
\begin{eqnarray}
t^{1/\alpha}|\XN\mu(t)|&\le& c\frac{k}{tH}\cdot \Bigl(\frac tk\Bigr)^ke^{k-t}t^{1/\alpha}\label{X9} \\
&\le & \e+ c_1 (k^{1/\alpha}/H)\cdot \chi_{\{t:k/2<t<2k\}},\qquad 0\le
t\le A,\nonumber
\end{eqnarray}
and
\begin{equation}
t^{1/\alpha}|\XN\mu(t)|\le \frac{c\tau e^{-t}}{H} \cdot
kt^{1/\alpha}e^{t/A} \le c_1\frac{A^{k-1+1/\alpha}\sqrt k e^{-A}}{H}\le
\e,\quad t\ge A. \label{X10}
\end{equation}

Now we fix $0<\psi<\beta-\frac \alpha 2$ and $k=\psi H^\alpha\log
H$ in such a way that $k\in\mathbb N$. Then \eqref{1} is satisfied
for large $H$, \eqref{X8} implies \eqref{X6}, and \eqref{X9},
\eqref{X10} imply \eqref{X5}.
\end{proof}

\begin{proof}[Proof of Theorem \rm\ref{example}]
Without loss of generality, we can assume that $\gamma$ is
non-increasing.

Our function $f$ will be defined by an inductive construction. Set
$Q_1=1$. On step $n\ge 1$ we use Lemma~\ref{prop} with $Q=Q_n$,
$\e=2^{-n}$ to find $\mu_n$, $k_n$ satisfying
\eqref{X1}--\eqref{X6}. Since
\begin{eqnarray*}
\lim_{|z|\to\infty,\, z\in \Omega}\XC\mu_n(z)&=&0,\\
\lim_{t\to\infty}\XL\mu_n(t)&=&0,\
\end{eqnarray*}
we can find $Q_{n+1}>k_n$ such that
\begin{eqnarray}
|\XC\mu_n(z)|&\le&\e,\qquad |z|>Q_{n+1},
\,z\in \Omega, \label{11}\\
|\XL\mu_n(t)|&\le&\e,\qquad t>Q_{n+1},
\end{eqnarray}
which completes the induction step.

Finally, for some numbers $\phi_n\in\mathbb C$, $|\phi_n|=1$, to be
chosen iteratively later on, we set
$$
f=\sum_{n\ge 1} \phi_n\gamma(k_n)^{1/\alpha}\XL \mu_n.
$$
The series above converges absolutely by the choice of $\mu_n,$ and
therefore  $f\in C_0(\mathbb R_+).$ Similarly, the function
$$
\widehat f=\sum_{n\ge 1} \phi_n\gamma(k_n)^{1/\alpha}\XC \mu_n
$$
extends analytically to $\Omega$, and
$$
\widehat f (z) (1 +|\Im z|)^{-\beta} \to 0, \qquad |z| \to
\infty, \,z\in \Omega.
$$\\
Finally, for $m\ge 1$ we have
\begin{eqnarray*}
&& \frac{k_m}{\gamma(k_m) \log k_m} \Bigl|\widehat
f(0)-\int_{0}^{k_m} f(s) \, ds \Bigr|^{\alpha}\\&=&
\frac{k_m}{\gamma(k_m) \log k_m} \Bigl| \sum_{n\ge 1}
\phi_n\gamma(k_n)^{1/\alpha}\XN \mu_n(k_m)
\Bigr|^{\alpha} \\
&\ge& \frac{c k_m}{\gamma(k_m) \log k_m} \Big( \Bigl| \sum_{1\le
n\le m} \phi_n\gamma(k_n)^{1/\alpha}\XN
\mu_n(k_m)\Bigr|^{\alpha}\\&-& \Bigl| \sum_{n>m}
\phi_n\gamma(k_n)^{1/\alpha}\XN \mu_n(k_m)\Bigr|^{\alpha} \Bigr).
\end{eqnarray*}

On step $m\ge 1$ we choose $\phi_m$ in such a way that
\begin{multline*}
\Bigl| \sum_{1\le n\le m} \phi_n\gamma(k_n)^{1/\alpha}\XN
\mu_n(k_m)\Bigr|\\= \Bigl|\gamma(k_m)^{1/\alpha}\XN \mu_m(k_m)
\Bigr| + \Bigl| \sum_{1\le n< m} \phi_n\gamma(k_n)^{1/\alpha}\XN
\mu_n(k_m)\Bigr|.
\end{multline*}
Then
\begin{gather*}
\frac{k_m}{\gamma(k_m) \log k_m} \Big( \Bigl| \sum_{1\le n\le m}
\phi_n\gamma(k_n)^{1/\alpha}\XN \mu_n(k_m)\Bigr|^{\alpha}\\- \Bigl|
\sum_{n>m} \phi_n\gamma(k_n)^{1/\alpha}\XN \mu_n(k_m)\Bigr|^{\alpha}
\Bigr)\\
\ge \frac{ k_m}{\gamma(k_m) \log k_m} \Big(\gamma(k_m)|\XN
\mu_m(k_m)|^{\alpha}- \Bigl| \sum_{n>m}
\phi_n\gamma(k_n)^{1/\alpha}\XN \mu_n(k_m)
\Bigr|^{\alpha}\Bigr)\\
\ge c_1-\frac{c_2}{\log k_m} \Bigl(\sum_{n>m}
2^{-n/\alpha}\Bigr)^\alpha \ge c_3>0.
\end{gather*}
\end{proof}

\begin{remark}
A variant of our construction works with $\gamma =1$ in Theorem
\ref{example}, if one looks just for $f \in L^{\infty}(\mathbb
R_+)$.
\end{remark}

\section{Decay of Banach space semigroups}\label{Banachh}

Using the construction of Theorem \ref{example}, we show next that
the analogue of Theorem \ref{battysem} for $C_0$-semigroups on
Banach spaces, Theorem \ref{battysem2}, is also sharp. Estimates
for the local resolvents of group generators similar to the ones
used below has been also employed, in particular, in \cite{ABHN01}
and \cite[Section II.4.6]{ChTo08}.

\begin{theorem}\label{xyz}
Given $\alpha > 0$, there exist a Banach space $X_{\alpha}$ and a
bounded $C_0$-semigroup $(T(t))_{t \ge 0}$ on $X_{\alpha}$ with
generator $A$  such that
$$
{\rm (a)}\qquad \|R(is,A)\|={\rm O}(|s|^{\alpha}), \qquad |s| \to \infty,
$$
and
$$
{\rm (b)}\qquad  \limsup_{t \to \infty} \Bigl(\frac{ t} { \log t}\Bigr)^{\frac{1}{\alpha}} \|T(t) A^{-1}\|
> 0.
$$
\end{theorem}

\begin{proof}
Let $(S(t))_{t \ge 0}$ be the left shift semigroup on ${\rm BUC}
(\mathbb R_+ )$, let $\Omega:=\{\lambda \in\mathbb C: \Re \lambda
> -1/(1+|\Im \lambda|)^\alpha\}$, and let $\Omega_0:=\Omega \cap
\{\lambda \in \mathbb C: |\Re \lambda| < 1 \}$. Furthermore, let
$X_{\alpha}$ be the space of functions $f\in {\rm BUC}(\R_+ )$
such that the Laplace transform $\widehat{f}$ extends to an
analytic function in $\Omega_0$, and
\begin{equation*}
|\widehat f(\lambda)|(1+|\Im \lambda|)^{-\alpha} \to 0, \qquad
\lambda \to \infty, \, \, \lambda \in \Omega_0.
\end{equation*}
Then $X_\alpha$ equipped with the norm
\begin{equation*}
\| f \|_{X_{\alpha}}  = : \| f \|_{\infty} + \|f \|_\alpha :=
\|f\|_{\infty} + \sup_{\lambda \in \Omega_0} |\widehat f(\lambda)|
(1+|\Im \lambda|)^{-\alpha},
\end{equation*}
is a Banach space. Moreover, $S(t)X_\alpha\subset X_\alpha$, $t
\ge 0$, and  the restriction $(T(t))_{t \ge 0}$ of $(S(t))_{t \ge
0}$ to $X_\alpha$ is also a $C_0$ -semigroup. To prove this assertion it
suffices to observe that
\begin{eqnarray*}
&& \widehat f(\lambda)-\widehat {T(t)f} (\lambda) \\
&=&\int_{0}^{\infty}e^{-\lambda s} f(s) \, ds - \int_{0}^{\infty}e^{-\lambda s} f (t+s) \, ds \\
&=& (1-e^{\lambda t})\widehat{f}(\lambda) + e^{\lambda t} \int_{0}^{t}
e^{-\lambda s}f(s) \, ds,\qquad \Re\lambda>0,
\end{eqnarray*}
and the same equality holds on $\Omega_0$. By the definition of $X_\alpha$,
$$
\| T(t)f - f \|_\alpha \to 0, \qquad t \to 0+,
$$
and then
$$
\| T(t)f - f \|_{X_\alpha} \to 0, \qquad t \to 0+.
$$

Let $A$ stand for
the generator of $(T(t))_{t \ge 0}$. 

Next we prove that for every $f \in X_\alpha$ the local resolvent
$R(\lambda,A)f$ satisfies the estimate
\begin{equation}\label{polyn}
\|R(\lambda,A)f\|_{X_\alpha}\le C\big(1+|\Im \lambda|
\big)^\alpha \| f\|_{X_\alpha}, \qquad 0<\Re\lambda<1.
\end{equation}
We will estimate the quantities  $\|R(\lambda,A)f\|_{\infty}$ and $\|R(\lambda,A)f\|_\alpha$
separately.

Observe first that for every $t\in\R_+$ and every $\lambda\in\C_+$
one has
$$
\bigl(R(\lambda,A)f \bigr)(t)= \int_0^\infty e^{-\lambda s} f(t+s) \; ds
= e^{\lambda t} \widehat{f} (\lambda ) - \int_0^t e^{\lambda (t-s)} f(s) \; ds .
$$
It follows that for every fixed $t\in\R_+$ the function
$\lambda\mapsto \left(R(\lambda, A)f \right) (t)$
extends to an analytic function on $\Omega_0$, and moreover
$$
|(R(\lambda,A)f)(t)| \le \left\{
\begin{array}{ll}
\frac{\| f\|_\infty}{|\Re \lambda|} & \text{if } \Re\lambda >0 , \\[2mm]
\frac{\| f\|_\infty}{|\Re \lambda|} + | \widehat{f}(\lambda)| &
\text{if } \Re\lambda <0 .
\end{array} \right.
$$

Applying Levinson's $\log$-$\log$ theorem (see, for example, \cite[VII D7]{Koo})
or, rather, its polynomial growth version \cite[Lemma 4.6.6]{ABHN01}
to $(R(\lambda,A)f)(t)$ in the squares $\{\lambda:|\Re \lambda|<(s+2)^{-\alpha},\,
|s-\Im \lambda|<(s+2)^{-\alpha}\}$,
we conclude that
\begin{equation}\label{esti}
\big \|R(\lambda,A)f \big \|_{\infty}  \le C(1+|\Im\lambda|)^{\alpha} \big( \| f\|_{\infty}+ \| f \|_\alpha \big),
\qquad  0<\Re\lambda<1.
\end{equation}

Fix $\lambda$ with $\Re\lambda \in (0,1)$. To estimate $\|R(\lambda,A)f\|_\alpha$ note that
\begin{eqnarray*}
\widehat{(R(\lambda,A)f )}(\mu) &=&\int_{0}^{\infty}
e^{-\mu t}
\int_{0}^{\infty}  e^{-\lambda s} f(t+s) \, ds \, dt\\
&=&\int_{0}^{\infty}e^{(\lambda -\mu)t}
\int_{t}^{\infty}e^{-\lambda s}
f(s) \, ds \,dt\\
&=& -\frac{1}{\lambda-\mu} \int_{0}^{\infty}e^{-\lambda s} f(s) \, ds
+\frac{1}{\lambda-\mu} \int_{0}^{\infty}e^{-\mu t} f(t) \, dt\\
&=& -\frac{\widehat f(\lambda) -\widehat f(\mu)}{\lambda - \mu}, \qquad \Re\mu>1.
\end{eqnarray*}
Therefore, $\widehat{R(\lambda,A)f}$ extends analytically to $\Omega_0$, and
$$
\widehat{(R(\lambda,A)f )}(\mu)=\begin{cases} -\frac{\widehat
f(\lambda) -\widehat f(\mu)}{\lambda - \mu}, \qquad \lambda \neq
\mu,\, \mu \in \Omega_0,\\
-{\widehat f\,}^\prime(\mu), \qquad \qquad \lambda=\mu.
\end{cases}
$$

Now, if
$|\lambda-\mu|\ge 1$, $\mu\in \Omega_0$,
then
$$
|\widehat{(R(\lambda,A)f)}(\mu)|\le|\widehat{f}(\lambda)|+\widehat{f}(\mu)|
\le c(1+|\Im \mu|)^{\alpha}(1+|\Im \lambda|)^{\alpha}\|f\|_\alpha.
$$
Furthermore, if
$$
1>|\lambda-\mu|\ge \frac{1}{2(1+|\Im \lambda|)^\alpha},\quad \mu\in \Omega_0,
$$
then we have
\begin{multline*}
|\widehat{(R(\lambda,A)f )}(\mu)|\le
2(1+|\Im \lambda|)^{\alpha} \bigl(|\widehat{f}(\lambda)|+|\widehat{f}(\mu)| \bigr)\\ \le
c(1+|\Im \lambda|)^\alpha(1+|\Im \mu|)^\alpha \|f\|_\alpha.
\end{multline*}
Finally, if
$$
|\lambda-\mu|\le\frac{1}{2(1+|\Im\lambda|)^\alpha},
$$
then, applying Cauchy's formula on the circle
$$
C_{\lambda}:=\{z\in \mathbb C:|z-\lambda|=\frac 23(1+|\Im \lambda|)^{-\alpha}\},
$$
we obtain that
$$
|\widehat{(R(\lambda,A)f)}(\mu)|\le
c(1+|\Im \mu|)^{\alpha} (1+|\Im \lambda|)^{\alpha}\|f\|_\alpha.
$$

Thus,
\begin{equation}\label{esto}
\big\|R(\lambda,A)f \big \|_\alpha  \le C(1+|\lambda|)^{\alpha}  \| f \|_\alpha, \qquad 0<\Re\lambda<1.
\end{equation}
The estimates \eqref{esti} and \eqref{esto} together give us
\eqref{polyn}. Since
$$
\big\|R(\lambda,A)\big\|_{X_{\alpha}}\ge \frac{1}{{\rm
dist}\,(\lambda,\sigma(A))}, 
$$
the estimate \eqref{polyn} implies that $i\mathbb R \subset \mathbb C\setminus \sigma (A)$,
and that
\begin{equation}\label{es}
\| R(\lambda,A)\|_{X_{\alpha}}\le C \big(1+|\Im \lambda|\big)^{\alpha}, \qquad \lambda \in i\mathbb R.
\end{equation}

Since $\sigma(A)\cap i\mathbb R=\emptyset$, by Theorem \ref{ingham} we obtain
$$
A^{-1}f = \lim_{t\to \infty} \big(A^{-1}f - T(t)A^{-1}f \big)=-\int_{0}^{\infty} T(t)f \, dt,
$$
for every $f \in X_\alpha$.

By Lemma~\ref{prop} there exist $f_n\in X_\alpha$, $f_n = \XL\mu_n$, and $k_n\to\infty$ as $n\to\infty$,
such that
$$
\|f_n\|_{X_\alpha}\le 1, \qquad |\XN\mu_n (k_n)| \ge C(\alpha)\Bigl(\frac{\log k_n}{k_n}\Bigr)^{1/\alpha},
\quad n \ge 1.
$$
Therefore,
\begin{gather*}
\|T(k_n)A^{-1}\|_{X_{\alpha}} \ge \|T(k_n )A^{-1} f_n\|_{X_{\alpha}}
\ge\Bigl\|\int_{0}^{\infty}f_n (\cdot + k_n+ r) \, dr\Bigr \|_{\infty}\\
\ge \Bigl|\int_{k_n}^{\infty} f_n (r) \, dr\Bigr|
\ge \Bigl|\widehat f_n(0)-\int_{0}^{k_n} f_n (s) \, ds \Bigr|
=  \Bigl|\XN\mu_{n} (k_n) \Bigr|\\
\ge C(\alpha) \left(\frac{\log k_n}{k_n}\right)^{1/\alpha}.
\end{gather*}
\end{proof}

\section{Acknowledgments}

We would like to thank C. J. K. Batty for sending us a preliminary
version of \cite{BaDu08} prior to its publication and A. M.
Gomilko for helpful discussions and remarks.

\end{document}